\documentclass[twocolumn]{autart}    
\usepackage{amsmath,amsfonts,amssymb}
\usepackage{graphicx}          

\numberwithin{equation}{section}

\begin{document}

\begin{frontmatter}

\title{Probability of constructing prediction model for observable of a dynamical process via time series\thanksref{footnoteinfo}} 

\thanks[footnoteinfo]{This paper was not presented at any IFAC 
meeting. Corresponding author Xiao-Song Yang.}

\author[HUST]{Xiao-Song Yang}\ead{yangxs@hust.edu.cn},    

\address[HUST]{School of Mathematics and Statistics
	Huazhong University of Science and Technology
	Wuhan, 430074, China}  

\begin{keyword}                           
Linear systems, Time series analysis, observable, prediction methods, dynamical systems   
\end{keyword}                             

\begin{abstract}                          
One fundamental problem in studying dynamical process is whether it is possible and how to construct prediction model for an unknown system via sampled time series, especially in the modern big data era. The research in this area is beneficial to experimentalists in physics,  chemistry, especially, in biological science, where it is hard to construct a prediction models by first principles. Therefore constructing prediction model for the observable of a complex system is of great practical significance in various areas of science and engineering. In the present paper, we show in terms of Lebesgue measure that, at least in the linear case, one can almost surely construct by linear algebra approach the prediction model about observable of unknown systems only via observed time series in the settings of discrete time and continuous time linear systems, even in the situation that one has no information about the underlying dynamical process and output function that generates the observed time series. Interestingly, we can obtain the characteristic polynomial of the unknown matrices for both of discrete time and continuous time linear systems with probability 1.
\end{abstract}

\end{frontmatter}


\section{Introduction}
In recent years, the topic on data driven learning of dynamical systems has received much attention from theoretical and engineering researchers \cite{S.Atkinson}\cite{S.L.Brunton}\cite{K.Champion}\cite{J.S.North}. The main theme is to construct (train) DNNs models for approximation of the evolution of state variables of the dynamical systems by virtue of the observed dataset. In dealing with learning dynamics of a system with unknown equation via DNNs, a typical problem is that the available data is just a subset of the state variables, or more generally, the observed output of the systems, which is usually a low dimensional valued function of the underlying state variables, as known as time series.

Since what we know is just the observable, which is an output of the system of interest, and its governing equation is unknown, it is natural to expect that one can predict the dynamical behavior of the system by its output variables, and this leads to the question of whether it is possible to construct a dynamic equation, i.e., prediction model for the observable (output variables) to describe its evolution. If this is the case, then theoretically one can study the dynamical behavior of the underlying system based on the dynamic equation for the observables (observable output). 

Along this line of research is the problem of reconstructing or predicting future behavior of an (unknown) dynamical system by the so-called time-delayed measurements of observables \cite{timedelay}\cite{Takens}\cite{meandimension}\cite{Y.Gutman}\cite{signal}, and will be called sampled time series in this paper as usual. This is one of the central topics in non-linear data analysis for unknown systems, which gives rise to challenging mathematical questions related to practical tools and algorithms. The current theme of this topic is to study possibility of constructing prediction model that is conjugate with the original (unknown) system just by observed time series, as called time-delayed measurements in the current physical literature.

While the aforementioned research focuses on the existence of prediction model for the  observable of the underlying unknown system, our interest in the present paper is how to obtain the prediction model about observed time variable in the situation where only observed time series is available, as well as likely-hood of dong this. This question is of significance, because only in the case that the properties of the prediction model are clearly investigated can one well predict dynamical behavior of the unknown original dynamical process. More importantly, short-term prediction is preferable in many practical situations. This renders constructing prediction model for observable of a dynamical process via time series a research area requiring more efforts.

Constructing prediction model for the output function of an unknown dynamical system via observed time series is of paramount importance from practical point of view. It is natural to expect that one can construct prediction model (which is in fact is a map) by some kind of interpolation approach, or an approximated model by various learning methods based on the available data as time series. For nonlinear process, the deep learning neural network seems to be a promising tool for this task. However, construction of prediction model for nonlinear process needs a great amount of data for training a deep neural network. 

As a first step, we will consider in this paper the above problem under the assumption that the underlying unknown system can be described by a linear autonomous systems. As for linear process, fortunately, we can resort to the theory of linear algebra and polynomials to provide a satisfactory framework for this problem. Now our goal is to study under what conditions and how often one can construct a dynamical model for the evolution of observed variable of underlying system that is unknown to the investigators in the case that only time series is available, and demonstrate that one can almost surely construct a prediction model about observed variable, i.e., the observable of unknown systems just by sampled time series in the settings of discrete time linear and continuous time linear systems. The results thus obtained would hopefully shed light on the determination of dynamic equations concerning evolution of outputs in more general settings.

\section{Observability and dynamic equation for the output variables in linear systems}
As shown in the sequel, whether a time series can effectively predict dynamics of an underlying system is deeply connected with the theory of observability of dynamical system. In particular, if the time series is generated by a single (scalar) output of a dynamical system, to fully reveal dynamical property of the system one needs the concept of observability of the system with the single output. For reader’s convenience, we just give a specific version of this concept, for general treatment of the observability, the reader can refer to \cite{FeedbackSys}.

Consider the following linear system with a scalar output.
\begin{equation}\label{eq2.1}
	x(i+1)=Ax(i),y_i=cx(i),x(i)\in \mathbb{R}^n,c^T \in \mathbb{R}^n
\end{equation}
Or more generally,
\begin{equation}\label{eq2.2}
	x(i+1)=Ax(i)+b, y_i=cx(i), x_i\in \mathbb{R}^n
\end{equation}
where $b\in \mathbb{R}^n$is a constant vector, and $y_i=cx(i)$is a (scalar or multivariable) output.

\begin{defn}\label{def2.1}
	A state $x\in \mathbb{R}^n$ is said to be observable, if there is a positive integer $k>0$ such that from the output $\{y_i:i=0,1,\cdots,k\}$ one can determine the state. If every state of (\ref{eq2.2}) is observable, then we say that (\ref{eq2.2}) is (completely) observable.
\end{defn}

The matrix
\begin{align*} \nonumber
	Q=
	\begin{pmatrix}
		c\\cA\\ \vdots \\cA^{n-1}
	\end{pmatrix}
\end{align*}
is called observability matrix. A classical theorem in control theory is that (\ref{eq2.2}) is observable if and only if the observability matrix has full rank, that is $rank(Q)=n$.

\begin{thm}\label{th2.2}
	If the observability matrix satisfies \\$rank(Q)=n$, then:
	\par
		$\mathrm{(i)}$ For (\ref{eq2.1}) we have the following dynamic equation for the output variable
		\begin{equation}
			y_n=G(y_0,y_1,\cdots, y_{n-1})^T, y_i\in \mathbb{R}, i=0,1,2,\cdots \nonumber
		\end{equation}
		where $G=cA^n Q^{-1}$
		
		$\mathrm{(ii)}$ For (\ref{eq2.2}) we have the following dynamic equation for the output variable
		\begin{equation}
			y_n=G(y_0,y_1,\cdots, y_{n-1})^T+\widetilde{b}, y_i\in \mathbb{R}, i=0,1,2,\cdots \nonumber
		\end{equation}
		where $G=cA^n Q^{-1}$ and 
		\begin{align*} \nonumber
			\widetilde{b}&=
			c(A^{n-1}+\cdots+E)b
			\\&-cA^nQ^{-1}
			\begin{bmatrix}
				0\\
				c\\
				c(A+E)\\
				\vdots\\ 
				c(A^{n-2}+A^{n-3}+\cdots+E)
			\end{bmatrix}
			b
		\end{align*}
\end{thm}

\textbf{Proof:}
     For every $x_0\in \mathbb{R}^n$, we have $y_0=cx_0, y_1=cx_1=cAx_0, \cdots y_{n-1}=cA^{n-1}x_0$. That is
	\begin{align*} \nonumber
		(y_0,y_1,\cdots, y_{n-1})^T=
		\begin{pmatrix}
			c\\cA\\ \vdots \\cA^{n-1}
		\end{pmatrix}
		x_0=Qx_0
	\end{align*}
	Therefore
	\begin{equation}
		x_0=Q^{-1}(y_0,y_1,\cdots,y_{n-1})^T \nonumber
	\end{equation}
	Since $y_n=cA^nx_0$, we have
	\begin{equation}
		y_n=cA^nx_0=cA^nQ^{-1}(y_0,y_1,\cdots,y_{n-1})^T\nonumber
	\end{equation}
	From the above equation we have the dynamical equation of output as 
	\begin{equation}\label{eq2.3}
		y_{m+n}=G(y_m, y_{m+1},\cdots,y_{m+n-1})^T,m\geq0
	\end{equation}
	and this is what we want for $\mathrm{(i)}$. One can prove $\mathrm{(ii)}$ in the same manner.
\hfill $\square$\par

In fact, we can have an explicit form for the dynamic equation (\ref{eq2.3}).
\begin{cor}\label{coro1}
	If the observability matrix satisfies $rank(Q)=n$, then the dynamical equation for the output variable is of the following form
	\begin{equation}\label{eq2.4}
		y_{m+n}=-a_0y_m - \cdots - a_{n-1}y_{m+n-1}, m\geq 0
	\end{equation}
\end{cor}

\textbf{Proof:}
	Let 
	\begin{equation}
		P(\lambda)=\lambda^n + a_{n-1}\lambda^{n-1} + \cdots + a_1\lambda+a_0 \nonumber 
	\end{equation}
	be the characteristic polynomial of the matrix, then by Cayley-Hamilton theorem, we have
	\begin{equation}
		P(A)=A^n + a_{n-1}A^{n-1} + \cdots + a_1A + a_0I \nonumber 
	\end{equation}
	It follows that
	\begin{equation}
		A^n = -a_{n-1}A^{n-1} - \cdots - a_1A - a_0I \nonumber 
	\end{equation}
	Now
	\begin{align*}
		G&=cA^nQ^{-1}\\
		 &=(-a_{n-1}cA^{n-1} - \cdots - a_1cA - a_0cI)Q^{-1}\\
		 &= (-a_0, -a_1, \cdots , -a_{n-1}) \nonumber
	\end{align*}
	We have
	\begin{align*}
		y_{m+n} &= G(y_m, y_{m+1}, \cdots , y_{m+n-1})^T \\
		&= -a_0y_m - a_1y_{m+1} - \cdots - a_{n-1}y_{m+n-1} 
	\end{align*}
\hfill $\square$\par

If we rewrite (\ref{eq2.3}) in a vector form as the following  n-dimensional linear system, then it is equivalent to (\ref{eq2.1}) in the sense of conjugation.

Let $z = (z_1,z_2, \cdots, z_n)^T$ with $z_1=c^T x, z_2=c^T Ax, \cdots$, and $z_n = c^T A^{n-1}x$, then the following fact is obvious.

\begin{prop}\label{prop1}
	System (\ref{eq2.1}) is conjugate with the following system
	\begin{equation}
		z(m+1) = A_cz(m) \nonumber
	\end{equation}
	where $A_c$ is the companion matrix of the following form
	\begin{align*}
		A_c=
		\begin{bmatrix}
			0 & 1 & 0 & \cdots & 0\\
			0 & 0 & 1 & \cdots & 0\\
			\vdots & \vdots &\vdots & \vdots & \vdots\\
			-a_0 & -a_1 & \cdots & \cdots & -a_{n-1}
		\end{bmatrix}
	\end{align*}
\end{prop}

Recall that two systems $x(i+1)=f(x(i))$ and $x(i+1)=g(x(i))$ are said to be conjugate, where $f,g:\mathbb{R}\rightarrow\mathbb{R}$ are continuous, if there is a homeomorphism $T:\mathbb{R}\rightarrow\mathbb{R}$ such that $f=T \circ g \circ T^{-1}$. Notice that two systems being conjugate means that they have qualitatively same dynamical behavior.

\begin{cor}
	If the time series generated by an observable output function satisfies the observability, then we can determine the stability of (\ref{eq2.1}) via prediction model for observed time series.
\end{cor}

\section{Constructing prediction model for observable via time series—deterministic case}
From practical viewpoint, it is more meaningful to have a prediction model for the observable of a (unknown) underlying system, i.e. the output of the underlying system.  In this section we consider for convenience the case of complete observable of (\ref{eq2.1}) in terms of an output. In this case the dynamic equation of the output is of the form of (\ref{eq2.4}), i.e.
\begin{equation}
	y_{m+n}=-a_0y_m - a_1y_{m+1} - \cdots - a_{n-1}y_{m+n-1}, m\geq0 \nonumber
\end{equation}
What we do is to estimate the values of $a_i, i=0,1,\cdots,n-1$ by observed time series.

In the deterministic setting, the observed time series is assumed to be without noise.

We have the following theorem.
\begin{thm}\label{th3.1}
	Assume that (\ref{eq2.2}) is observable with the output, then a time series $\{y_j=cA^jx_0:j=0,1,\cdots\}$ generated by an output $y_i=cx(i)$ of (\ref{eq2.1}), can determine the value $a_i, i=0,1,\cdots,n-1$ provided $\vert A \vert \neq 0$, and the vectors $A^jx_0, j=0,1,\cdots,n-1$, are linear independent.
\end{thm}

\textbf{Proof:}
	For the given time series $\{y_j=cA^jx_0:j=0,1,\cdots\}$, consider the linear equations
	\begin{align*} 
		y_n&=-a_0y_0-a_1y_1-\cdots-a_{n-1}y_{n-1} \nonumber\\
		y_n+1&=-a_0y_1-a_1y_2-\cdots-a_{n-1}y_n \nonumber\\
		\vdots& \nonumber\\
		y_{2n-1}&=-a_0y_{n-1}-a_1y_n-\cdots-a_{n-1}y_{2n-2}  \nonumber
	\end{align*}
	Let
	\begin{align*} 
		H=
		\begin{bmatrix}
			y_0 & y_1 & \cdots & y_{n-1} \nonumber\\
			y_1 & y_2 & \cdots & y_n \nonumber\\
			\vdots & \vdots & \vdots & \vdots \nonumber\\
			y_{n-1} & y_n &\cdots &y_{2n-2} \nonumber
		\end{bmatrix}
	\end{align*}
	We have
	\begin{align*}
		H&=
		\begin{bmatrix}
			cx_0 & cAx_0 & \cdots & cA^{n-1}x_0\\
			cAx_0 & cA^2x_0 & \cdots & cA^nx_0\\
			\vdots & \vdots & \vdots & \vdots\\
			cA^{n-1}x_0 & cA^nx_0 & \cdots & cA^{2n-2}x_0
		\end{bmatrix}\\
		&=
		\begin{bmatrix}
			c\\
			cA\\
			\vdots\\
			cA^{n-1}
		\end{bmatrix}
		\begin{bmatrix}
			x_0 & Ax_0 & \cdots & A^{n-1}x_0
		\end{bmatrix} \\
		&=QM(A,x_0) \nonumber
	\end{align*}
	Then the conditions assumed in the theorem imply that $\vert Q \vert \neq 0$ and $\vert M(A,x_0) \vert \neq 0$. Therefore we can determine each $a_i$ by virtue of $\{y_j:j=0,1,\cdots,2n-2\}$.
	
	Now for any integer $k > 0$ consider the following equations 
	\begin{align*} 
		y_{n+k}&=-a_0y_k-a_1y_{k+1}-\cdots-a_{n-1}y_{k+n-1} \nonumber\\
		y_{n+k+1}&=-a_0y_{k+1}-a_1y_{k+2}-\cdots-a_{n-1}y_{k+n} \nonumber\\
		\vdots& \nonumber\\
		y_{2n+k-1}&=-a_0y_{k+n-1}-a_1y_{k+n}-\cdots-a_{n-1}y_{k+2n-2} \nonumber
	\end{align*}
	Let
	\begin{align*} \nonumber
		H_k=
		\begin{bmatrix}
			y_k & y_{k+1} & \cdots & y_{k+n-1}\\
			y_{k+1} & y_{k+2} & \cdots & y_{k+n}\\
			\vdots&\\
			y_{n+k-1} & y_{n+k} & \cdots & y_{2n+k-2}
		\end{bmatrix}
	\end{align*}
	Then $\vert H_k \vert \neq 0$, thus we have
	\begin{align*} \nonumber
		\begin{pmatrix}
			-a_0\\
			-a_1\\
			\vdots\\
			-a_{n-1}
		\end{pmatrix}
		=H_k^{-1}
		\begin{pmatrix}
			y_{n+k}\\
			y_{n+k+1}\\
			\vdots\\
			y_{2n+k-1}
		\end{pmatrix}
	\end{align*}
	It is easy to see that
	\begin{align*}
		H_k&=
		\begin{bmatrix}
			c\\
			cA\\
			\vdots\\
			cA^{n-1}
		\end{bmatrix}
		A^k
		\begin{bmatrix}
			x_0 & Ax_0 & \cdots & A^{n-1}x_0
		\end{bmatrix} \\
		&=QA^kM(A,x_0) \nonumber
	\end{align*}
	In view of the assumptions of the above theorem, the matrix $H_k$ constructed by $\{y_j:j=k+1, \cdots, 2n+k-2\}$ is invertible, thus one can determine each $a_i$ by any piece $\{y_j:j=k+1, \cdots, 2n+k-2\}$ of the time series $\{y_j=cA^j x_0: j=0, 1, \cdots\}$.
\hfill $\square$\par

\begin{rem}\label{remark1}
	A question thus arises that if it is always possible that $\{y_j=cA^j x_0: j=0, 1, \cdots\}$ satisfies $\vert M(A,x_0) \vert \neq 0$ for $x_0 \neq 0$. Clearly this is not always the case, for if $x_0 \neq 0$ is a real eigenvector of $A$, then $\vert M(A,x_0) \vert=0$. Fortunately, it can be easily proved that the set of points in $\mathbb{R}^n$ satisfying $\vert M(A,x_0) \vert=0$ is open and dense in $\mathbb{R}^n$.  In the next section, we will prove that such a set is of full Lebesgue measure.
\end{rem}

\section{Probability of detecting dynamics of discrete time linear systems}

First we give the following definition for convenience of the sequel arguments

\begin{defn}\label{def4.1}
	Let $D\subset \mathbb{R}^n$ be a set with nonempty interior. A set $B \subset D$ is said to be of full (Lebesgue) measure with respect to $D$, if $D-B$ has zero Lebesgue measure.
\end{defn}
The following fact is obvious.
\begin{lem}\label{lemma4.2}
	Let $M_{n \times n}$ be the matrix space of $n \times n$,  which can be equivalently regarded as the Euclidean space $\mathbb{R}^{n \times n}$. Let $\bar{M}_{n \times n}$ be the subset of the matrices that have distinct eigenvalues. Then $\bar{M}_{n \times n}$ is dense and open in $M_{n \times n}$.
\end{lem}

It is not satisfactory from practical point of view that a property holds in an open and dense subset of a matric space, because the openness and density do not necessarily guarantee full measure, and such a set may have very small measure on a set of interest, therefore it is of practical significance to study a property holding on a set with full measure.  Thus the following observation is more meaningful.

\begin{lem}\label{lemma4.3}
	$\bar{M}_{n \times n}$ has full Lebesgue measure in $M_{n \times n}$.
\end{lem}

\textbf{Proof:}
	Consider matrix $A\in M_{n \times n}$, and its characteristic polynomial
	\begin{equation}
		P_{A}(\lambda)=\lambda^n+a_1\lambda^{n-1}+\cdots+a_{n-1}\lambda+a_n \nonumber
	\end{equation}
	Then 
	\begin{equation}
		{P_{A}^{'}}(\lambda)=n\lambda^{n-1}+(n-1)a_1\lambda^{n-2}+\cdots+a_{n-1} \nonumber
	\end{equation}
	The discriminant $D(P_A)$ of $P_A$ is 
	\begin{equation}
		D(P_{A})=(-1)^{\frac{n(n-1)}{2}}R\left(P_{A}, P_{A}^{'} \right) \nonumber
	\end{equation}
	where $R\left(P_{A}, P_{A}^{'}\right)$ is the so-called resultant of $P_{A}$ and $P_{A}^{'}$ \cite{usingAG}, which is the determinant of an $\left(n + (n-1)\right) \times \left(n + (n-1)\right)$ matrix with all its entries are the coefficients of $P_{A}$, thus is a polynomial function of the coefficients of $P_{A}$. Notice that every coefficient $a_i$
	is also a polynomial function of some entries of $A$, it is easy to see that $R\left(P_{A}, P_{A}^{'}\right)$ is an analytic function of the entries of $A$. Therefore $D(P_{A}):M_{n \times n}\rightarrow \mathbb{R}$ is an analytic function. 
	
	Note that the function $D(P_{A})$ is not zero if and only if $A$ has $n$ distinct eigenvalues and the set of zero points of $D(P_{A})$ is a union of submanifolds with dimensions lower than $n \times n$. Thus the set of zero points of $D(P_{A})$ has zero Lebesgue measure in the space $M_{n \times n}$ and this implies that $\bar{M}_{n \times n}$ has full Lebesgue measure in $M_{n \times n}$ in the sense of \textbf{Definition \ref{def4.1}}.
\hfill $\square$\par

\begin{rem}\label{remark2}
	Let $M_c \subset M_{n \times n}$ be the set of matrices that are similar to companion matrices, then it is easy to see in view of elementary matrix theory that $M_c$ is open and dense in $M_{n \times n}$ amd $\bar{M}_{n \times n}\subset M_{c}$. Furthermore  the above lemma implies that $M_c$ has full Lebesgue measure in $M_{n \times n}$
\end{rem}

\begin{lem}\label{lemma4.4}
	If $A\in \bar{M}_{n \times n}$, then there is a vector $v \in \mathbb{R}^n$ such that $v, Av, \cdots, A^{n-1}v$ are linearly independent; and there is a vector $c\in \mathbb{R}^n$ such that $c^T, A^T c^T, \cdots, \left(A^{n-1}\right)^T c^T$ are linearly independent.
\end{lem}

\textbf{Proof:}
	Since a matrix $A\in \bar{M}_{n \times n}$ has $n$ different eigenvalues, it is similar to a companion matrix $H$, i.e., there exists a nonsingular matrix $T\in M_{n \times n}$ such that $TAT^{-1}=H$. For $H$, there is a $x\in \mathbb{R}^n$ such that $x, Hx, \cdots, H^{n-1}x$ are linearly independent. This implies that the matrix $[x$, $TAT^{-1}x$, $\cdots$, $\left( TAT^{-1}\right)^{n-1}x]=$ $T[T^{-1}x$, $AT^{-1}x$,$ \cdots$, $A^{n-1}T^{-1}x]$ has full rank, therefore $T^{-1}x$, $AT^{-1}x$,$ \cdots$, $A^{n-1}T^{-1}x$ are linearly independent. Let $v=T^{-1}x$, we complete the proof. The second statement can be proved in a similar manner. 
\hfill $\square$\par

\begin{defn}\label{def4.5}
	Consider a matrix $A\in M_{n \times n}$. If a vector $v\in \mathbb{R}^n$ satisfies that $v, Av, \cdots. A^{n-1}v$ are linearly independent, then the vector $v$ is said to have $A$-linear independence property. If a vector $c\in \mathbb{R}^n$ satisfies that $c^T, A^T c^T, \cdots, \left(A^{n-1}\right)^T c^T$ are linearly independent. Then the vector $c$ is said to have $A^T$-linear independence property.
\end{defn}

\begin{rem}\label{remark3}
	A matrix having the above property is also called a companion matrix. In control theory, $c^T, A^T c^T, \cdots, \left(A^{n-1}\right)^T c^T$ being linearly independent is equivalent to complete observability of system (\ref{eq2.1}).
\end{rem}

\begin{prop}\label{prop2}
	For each $A\in \bar{M}_{n \times n}$, the set of vectors that satisfies $A$-linear independence property is open and dense in $\mathbb{R}^n$, and set of vectors that satisfying $A^T$-linear independence property is also open and dense in $\mathbb{R}^n$.
\end{prop}

\textbf{Proof:}
	The openness is obvious. To prove the denseness, suppose on the contrary, that there is an open subset $U$ of $\mathbb{R}^n$ such that the $\det [x, Ax, \cdots, A^{n-1}x]=0$, for $\forall x\in U$. Since $\det [x, Ax, \cdots, A^{n-1}x]=0$ is an analytic function, one has $\det [x, Ax, \cdots, A^{n-1}x]\equiv 0$, for $\forall x\in \mathbb{R}^n$. On the other hand, $A\in \bar{M}_{n \times n}$ implies that it has $n$ distinct eigenvalues, and there is a vector $v\in \mathbb{R}^n$ such that $\det [v, Av, \cdots, A^{n-1}v]\neq 0$, leading to a contradiction. The same argument can apply to the case of $A^T$-linear independence property.
\hfill $\square$\par

In the following we now give a fact on  $A$-linear independence property in terms of probability.

\begin{prop}\label{prop3}
	For any compact domain $K\subset \mathbb{R}^n$ and a matrix $A\in \bar{M}_{n \times n}$, the set of points $x\in K$ satisfying $A$-linear independence property is of full Lebesgue measure with respect to $K$. Therefore the Lebesgue-induce probability that $x\in K$ satisfies $A$-linear independence property is 1.
\end{prop}

\textbf{Proof:}
	Note that if $f$ is an analytic function defined on an $n$-dimensional domain, then its zero set $Z_f$ can be written as the union lower dimensional (real analytic) submanifolds \cite{AnaPoly}, thus the Lebesgue measure of $Z_f$ is zero in $\mathbb{R}^n$. This implies that $K-Z_f$ has full Lebesgue measure on $K$.
\hfill $\square$\par

Similarly we have the following fact.

\begin{prop}\label{prop4}
	For any Lebesgue measurable set (compact domain) $K\subset \mathbb{R}^n$ and a matrix $A\in \bar{M}_{n \times n}$, the Lebesgue-induced probability that a point $c\in K$ satisfying $A^T$-linear independence property is one.
\end{prop}

\begin{rem}\label{remark4}
	This fact combined with Lemma \ref{lemma4.4}, implies that randomly choosing $c\in K$ and $A\in M_{n \times n}$, the system $(A, c)$ is completely observable with full probability, i.e., $(A, c)$ is almostly surely completely observable.
\end{rem}

Now summarizing the above facts we have the following main theorem.

\begin{thm}\label{th4.6}
	For any Lebesgue measurable set (compact domain) $K\subset \mathbb{R}^n \times M_{n\times n} \times \mathbb{R}^n$, the Lebesgue induce probability that the dynamic of a system can be detected by time series generated by $(c, A, x)\in K$ is 1.
\end{thm}

\textbf{Proof:}
	It is enough to consider Lebesgue measurable sets $K_1 \subset \mathbb{R}^n, K_2\subset M_{n\times n}$ and $K_3\subset \mathbb{R}^n$. In view of lemma \ref{lemma4.3}, the subset $\bar{M}_{n \times n}$ of the matrices that have $n$ distinct eigenvalues is dense and open in $M_{n \times n}$ and has full Lebesgue measure in $M_{n \times n}$.
\hfill $\square$\par

Therefore, one can almost surely construct a dynamical model about the observed output of an unknown dynamical system just via observed output time series, and the constructed model is the same as the original unknown system in the sense of the conjugation!

\section{Detecting dynamics of continuous time linear systems}
Consider the following linear system with a scalar output.
\begin{equation}\label{eq5.1}
	\dot{x}=Ax, y=cx, x\in \mathbb{R}^n
\end{equation}
It is well known that (\ref{eq2.1}) is observable if and only if the so-called observability matrix.
\begin{align*} \nonumber
	Q=
	\begin{pmatrix}
		c\\
		cA\\
		\vdots\\
		cA^{n-1}
	\end{pmatrix}
\end{align*}
has full rank, that is $rank(Q)=n$. Now we know that for almost $c\in \mathbb{R}^n$ and $A\in M_{n \times n}$, the observability matrix has full rank. It is of practical significance to detect dynamics (\ref{eq5.1}) via sampled time series $y(i)=c\phi (i\lambda, x_0)$, where $\phi(t, x_0)$ is the solution of (\ref{eq2.1}) with initial state $x_0$ and $\lambda$ is a small positive number.

Note that $\phi(t, x_0)=e^{tA}x_0$, from (\ref{eq5.1}) we have a discrete time system
\begin{equation}\label{eq5.2}
	x(i+1)=Bx(i), y_i=cx(i), x(i)\in \mathbb{R}^n
\end{equation}
with $B=e^{\lambda A}$.It is an approximation of (\ref{eq5.1}).

Let
\begin{align*} \nonumber
	P=
	\begin{pmatrix}
		c\\
		cB\\
		\vdots\\
		cB^{n-1}
	\end{pmatrix}
\end{align*}
be the observability matrix of (\ref{eq2.2}). Now we show that following fact.

\begin{prop}\label{prop5}
	The two observability matrices satisfy $rank(Q)=rank(P)$ if $\lambda$ is sufficiently small.
\end{prop}

\textbf{Proof:}
	Consider the following $n$th Taylor formula
	\begin{align*} 
		e^{\lambda A}&=I + \lambda A +\frac{\lambda^2 A^2}{2!}+\cdots+\frac{\lambda^{n-1} A^{n-1}}{(n-1)!}\\
		&+o(\lambda^{n-1}) \nonumber\\
		e^{2\lambda A}&=I + 2\lambda A +\frac{(2\lambda)^2 A^2}{2!}+\cdots+\frac{(2\lambda)^{n-1} A^{n-1}}{(n-1)!}\\
		&+o(\lambda^{n-1}) \nonumber\\
		\vdots& \nonumber\\
		e^{(n-1)\lambda A}&=I + (n-1)\lambda A +\frac{((n-1)\lambda)^2 A^2}{2!}+\cdots\\
		&+\frac{((n-1)\lambda)^{n-1} A^{n-1}}{(n-1)!}+o(\lambda^{n-1}) \nonumber
	\end{align*}
	Then it is easy to see that 
	\begin{align*}
		P&=
		\begin{pmatrix}
			c\\
			cB\\
			\vdots\\
			cB^{n-1}
		\end{pmatrix}
		=
		\begin{pmatrix}
			c\\
			ce^{\lambda A}\\
			\vdots\\
			ce^{(n-1)\lambda A} \\
		\end{pmatrix}  \nonumber\\
		&=\begin{bmatrix}
			1 & 0 & \cdots & 0\\
			1 & \lambda & \cdots & \frac{\lambda^{n-1}}{(n-1)!}\\
			\vdots & \vdots & \vdots & \vdots\\
			1 & (n-1)\lambda & \cdots & \frac{((n-1)\lambda)^{n-1}}{(n-1)!} \\
		\end{bmatrix}
		\begin{pmatrix}
			c\\
			cA\\
			\vdots\\
			cA^{n-1}
		\end{pmatrix}
		+ o(\lambda^{n-1})  \nonumber\\
		&=T
		\begin{pmatrix}
			c\\
			cA\\
			\vdots\\
			cA^{n-1}
		\end{pmatrix}
		+o(\lambda^{n-1})  \nonumber
	\end{align*}
	It is obvious that $\det (T)\neq 0$. To see this, note that
	\begin{align*} 
		\det (T)&=
		\begin{vmatrix}
			\lambda & \frac{\lambda^2}{2!} & \cdots & \frac{\lambda^{n-1}}{(n-1)!}\\
			2\lambda & \frac{2^2 \lambda^2}{2!} & \cdots & \frac{2^{n-1} \lambda^{n-1}}{(n-1)!}\\
			\vdots & \vdots & \vdots & \vdots\\
			(n-1)\lambda & \frac{((n-1)\lambda)^2}{2!} & \cdots & \frac{((n-1)\lambda)^{n-1}}{(n-1)!}
		\end{vmatrix} \\
		&=\frac{1}{2!\cdots(n-1)!}\times\\
		&
		\begin{vmatrix}
			\lambda & \lambda^2 & \cdots & \lambda^{n-1}\\3
			2\lambda & 2^2\lambda^2 & \cdots & 2^{n-1}\lambda^{n-1}\\
			\vdots & \vdots & \vdots & \vdots \\
			(n-1)\lambda & ((n-1)\lambda)^2 & \cdots & ((n-1)\lambda)^{n-1}
		\end{vmatrix} \\
		&=\frac{1}{2!\cdots(n-1)!}\times H
	\end{align*}
	
	Note that
	\begin{eqnarray} 
		H&=
		\begin{vmatrix}
			\lambda & \lambda^2 & \cdots & \lambda^{n-1}\\
			2\lambda & 2^2\lambda^2 & \cdots & 2^{n-1}\lambda^{n-1}\\
			\vdots & \vdots & \vdots & \vdots \\
			(n-1)\lambda & ((n-1)\lambda)^2 & \cdots & ((n-1)\lambda)^{n-1}
		\end{vmatrix} \nonumber\\
		&=(n-1)!\lambda
		\begin{vmatrix}
			1 & \lambda & \cdots & \lambda^{n-2}\\
			1 & 2\lambda & \cdots & (2\lambda)^{n-2}\\
			\vdots & \vdots & \vdots & \vdots \\
			1 & (n-1)\lambda & \cdots & ((n-1)\lambda)^{n-2} \nonumber
		\end{vmatrix}
	\end{eqnarray}
	
	In view of the property of Vandermonde determinant, we have $H \neq 0$, since $\lambda > 0$. It follows that if $\lambda >0$ is small enough, then $P=[c^T, B^T c^T, \cdots, (B^{n-1})^T c^T]$ is nonsingular.
\hfill $\square$\par

It can be seen from the proof of the above theorem that if $A$ is similar to a companion matrix, then for sufficiently small $\lambda>0$, the matrix $e^{\lambda A}$ is also similar to some companion matrix, thus we have the following fact about continuous time linear system. Let $\phi (t,x)$ be the solution of (\ref{eq5.1}) with initial point $x\in \mathbb{R}^n$, this fact can be stated as follows.

\begin{thm}\label{th5.1}
	If $(A, c)$ of system (\ref{eq5.1}) is completely observable, then for $\lambda$ small enough, and one can construct a prediction model about time series ${y_i}$ generated from the output $y(t)$ via the sampled-time series $y_i=ce^{i \lambda A}x_0, i=0, 1, \cdots$, for almost of $x_0\in \mathbb{R}^n$, such that this model is conjugated with $x(i+1)=\phi(\lambda, x(i))$ generated (\ref{eq5.1}) in the sense of \textbf{Proposition \ref{prop1}}, as long as $\lambda>0$ is sufficiently small.
\end{thm}

\textbf{Proof:}
	It is easy to see that the  matrix $[c^T$, $B^T c^T$,$\cdots$, $(B^{n-1})^T c^T]^T$ with $B=e^{\lambda A}$ is nonsingular, showing that $P$ is similar to some companion matrix. From the arguments above, it can be asserted that for almost $x_0\in \mathbb{R}^n$ we can construct a prediction model for the output $y(t)$ by means of the sampled time series $y_i=ce^{i\lambda A}x_0, i=0,1,\cdots$. It is obvious that we can also solve the characteristic polynomial equation of the matrix $B=e^{\lambda A}$ via time series $\{y_i\}$ in view of Theorem \ref{th2.2}.
\hfill $\square$\par

This theorem implies that we can predict future behavior of output $y(t)$ of a linear process described by ordinary differential equations without resort to information of the ordinary differential equations. Also we can obtain the information of characteristic polynomial of the unknown  matrix $A$ by virtue of the relation $B=e^{\lambda A}$.

As a consequence of the above theorem and previous arguments, we have a general conclusion as follows.

\begin{thm}\label{th5.2}
	In linear case with the dimension of state space known in advance, even for an unknown continuous time system one can construct a continuous time dynamical system for the observed output of the original unknown system by discrete time series sampled from the output. And this is almost surely possible even we do not know the exact form of the output function.
\end{thm}

\section{Discussion}
In light of the discussions above, an elegant fact about linear dynamical process with  interactive state parameters, is that it is easy to construct by means of linear algebraic techniques a prediction model for the observable of this process. Furthermore, a consecutive time series of length not greater than $2n$ is enough to finish this task with full Lebesgue measure, or probability 1 in a suitable measure.

Although it is not the case that the above fact holds for general nonlinear process, the main result is meaningful to bring more hope to deep investigation of complex process via time series from probability view point.

\begin{ack}                               
\end{ack}

\bibliographystyle{plain}        
\bibliography{ReferenceFile}     

\begin{thebibliography}{10}

\bibitem{FeedbackSys}
Karl~Johan Astrom and Richard M.Murray.
\newblock {\em Feedback Systems: An Introduction for Scientists and Engineers,
  Second Edition}.
\newblock Princeton University Press, USA, 2021.

\bibitem{S.Atkinson}
Steven Atkinson, Waad Subber, Liping Wang, Genghis Khan, Philippe Hawi, and
  Roger Ghanem.
\newblock Data-driven discovery of free-form governing differential equations.
\newblock {\em arXiv e-prints}, pages 1--7, 2019.

\bibitem{timedelay}
Krzysztof Barański, Yonatan Gutman, and Adam Śpiewak.
\newblock On the shroer–sauer–ott–yorke predictability conjecture for
  time-delay embeddings.
\newblock {\em Communications in Mathematical Physics}, 391(2):609--641,
  February 2022.

\bibitem{S.L.Brunton}
Steven~L. Brunton, Joshua~L. Proctor, and J.~Nathan Kutz.
\newblock Discovering governing equations from data by sparse identification of
  nonlinear dynamical systems.
\newblock {\em Proceedings of the National Academy of Sciences},
  113(15):3932--3937, March 2016.

\bibitem{K.Champion}
Kathleen Champion, Bethany Lusch, J.~Nathan Kutz, and Steven~L. Brunton.
\newblock Data-driven discovery of coordinates and governing equations.
\newblock {\em Proceedings of the National Academy of Sciences},
  116(45):22445--22451, 2019.

\bibitem{usingAG}
David~A. Cox, John Little, and Donal O’shea.
\newblock {\em Using Algebraic Geometry}.
\newblock Springer New York, NY, USA, 2005.

\bibitem{Takens}
Yonatan Gutman.
\newblock Takens' embedding theorem with a continuous observable.
\newblock {\em arXiv e-prints}, pages 134--141, 2016.

\bibitem{meandimension}
Yonatan Gutman and Lei Jin.
\newblock Mean dimension and an embedding theorem for real flows.
\newblock {\em arXiv e-prints}, pages 161--181, 2020.

\bibitem{Y.Gutman}
Yonatan Gutman, Yixiao Qiao, and Gábor Szabó.
\newblock The embedding problem in topological dynamics and takens’ theorem.
\newblock {\em Nonlinearity}, 31(2):597--620, January 2018.

\bibitem{signal}
Yonatan Gutman, Yixiao Qiao, and Masaki Tsukamoto.
\newblock Application of signal analysis to the embedding problem of
  $\mathbb{Z}^k$-actions.
\newblock {\em arXiv e-prints}, 2017.

\bibitem{J.S.North}
Joshua~S. North, Christopher~K. Wikle, and Erin~M. Schliep.
\newblock A review of data-driven discovery for dynamic systems.
\newblock {\em arXiv e-prints}, 2022.

\bibitem{AnaPoly}
Q~I Rahman and G~Schmeisser.
\newblock {\em {Analytic Theory of Polynomials}}.
\newblock Oxford University Press, 09 2002.

\end{thebibliography}



\end{document}